\documentclass{amsart}
\usepackage[matrix, arrow, curve]{xy}
\usepackage{verbatim}
\usepackage{amsmath,amsfonts,amssymb,amsthm,epsfig}

\newcommand{\Crystals}{\mathsf{Crystals}}
\newcommand{\Cat}{\mathcal{C}}
\newcommand{\GLn}{\mathfrak{gl}_n}

\newcommand{\wt}{\operatorname{wt}}
\newcommand{\Hom}{\operatorname{Hom}}

\newcommand{\flip}{\operatorname{flip}}
\newcommand{\qg}{U_q(\mathfrak{g})}
\newcommand{\R}{{\mathbb R}}

\newtheorem{Theorem}{Theorem} 
\newtheorem{Proposition}{Proposition} 
\newtheorem{Lemma}{Lemma}

\newtheorem{Example}{Example}

\begin{document}

\title[Crystals and Coboundary Categories]{Crystals and Coboundary Categories}

\author{Andr\'e Henriques}
\email{andrhenr@math.mit.edu}
\address{Department of Mathematics\\ MIT \\ Cambridge, MA}

\author{Joel Kamnitzer}
\email{jkamnitz@math.berkeley.edu}
\address{Department of Mathematics\\ UC Berkeley \\ Berkeley, CA}

\begin{abstract}
Following an idea of A. Berenstein, we define a commutor for the category of crystals of a finite dimensional complex reductive Lie algebra.  We show that this endows the category of crystals with the structure of a coboundary category.  
Similar to the role of the braid group in braided categories, a group naturally acts on multiple tensor products in coboundary categories.  We call this group the cactus group and identify it as the fundamental group of the moduli space of marked real genus zero stable curves.
\end{abstract}

\date{\today}
\maketitle
\tableofcontents
\newpage

\section{Introduction}

\subsection{A commutor for crystals}
Let $ \mathfrak{g} $ be a finite dimensional complex reductive Lie algebra.  Crystals were introduced by Kashiwara as a combinatorial structure arising from the $ q \rightarrow 0 $ limit of a representation of the quantum group $ U_q(\mathfrak{g}) $.  Roughly, a crystal is a directed graph where the edges are labelled by simple roots of $ \mathfrak{g} $ and the vertices are labelled by weights of $ \mathfrak{g}$.  Crystals can be tensored together to produce a crystal whose underlying set is the product of the two underlying sets.  This tensor product is not symmetric in the sense that the map $ (a,b) \mapsto (b,a) $ is not an isomorphism from $A\otimes B$ to $B\otimes A$.

In this work, following an idea of Arkady Berenstein we construct natural isomorphisms $ \sigma_{A,B}: A \otimes B \rightarrow B \otimes A $ for crystals of finite-dimensional reductive Lie algebras.  The basic idea is to first produce a involution $ \xi: B \rightarrow B $ for each crystal $B$, which flips the crystal by exchanging highest weight elements with lowest weight elements.  Then we define: 
\begin{equation*}
\sigma(a,b) = \xi \big( \xi(b), \xi(a) \big) 
\end{equation*}
We call this natural isomorphism $\sigma$ the \textbf{commutor}. 

The commutor which we construct for crystals does not obey the axioms for a braided monoidal category. Instead we see that:
\begin{enumerate} 
\item $ \sigma_{B,A} \circ \sigma_{A,B} = 1 $
\item The following diagram commutes:
\begin{equation*} 
\xymatrix{
A \otimes B \otimes C \ar[d]_{\sigma_{A,B} \otimes 1} \ar[r]^{1 \otimes \sigma_{ B, C}} & A \otimes C \otimes B \ar[d]^{\sigma_{A, C \otimes B}}\\
B \otimes A \otimes C \ar[r]_{\sigma_{B \otimes A, C}} & C \otimes B \otimes A \\
}
\end{equation*}
\end{enumerate}

Our method requires $ \mathfrak{g} $ to be finite-dimensional since we use the long element of the Weyl group in order to produce the involution $ \xi $.  It would be interesting to see if there exists a commutor satisfying (i), (ii) for the category of highest weight crystals of an affine (or more generally Kac-Moody) Lie algebra.

\subsection{Relation with quantum groups}

Recall that the category of representations of $ \qg$ has a natural braiding which is constructed using the universal R matrix.  This braiding is not symmetric but it does obey the hexagon axiom and so endows the category of representations of $ \qg $ with the structure of a braided monoidal category.  

So, it is perhaps surprising that the commutor that we constructed for crystals is not braided.  However, following an analogous procedure to the one used for crystals, we construct a new commutor for the representations of $ \qg $ which obeys (i), (ii).  Using the canonical basis, we relate this commutor to our commutor for crystals.

\subsection{Coboundary categories}

In \cite{Drin}, Drinfel'd defined a coboundary category to be a monoidal category along with a commutor satisfying (i), (ii) above.  This name was chosen because the representation categories of coboundary Hopf algebras form coboundary categories. 

Interestingly, we stumbled across the axioms of a coboundary category in a different context: the hives and octahedron recurrence of Knutson, Tao and Woodward \cite{ktw}.  See \cite{us2} for the relation between these axioms and certain properties of the octahedron recurrence.  The connection to crystals was made when we found an equivalence of categories between the category of hives and the category of $\GLn$ crystals (see \cite{us}). 

Motivated by the example of crystals, we study the structure of coboundary categories. In braided cactgories, we have the braid group which acts naturally on multiple tensor products. There is an analogue of it for coboundary categories. We identify its structure, and call it the \bf cactus group\rm. It is the group which acts naturally on multiple tensor products of objects in a coboundary category.

In a coboundary category, the basic maps between repeated tensor products are reversals of intervals:
\begin{equation*}
s_{p,q} : A_1 \otimes \cdots \otimes A_n \rightarrow A_1 \otimes \cdots \otimes A_{p-1} \otimes A_q \otimes A_{q-1} \otimes \cdots \otimes A_p \otimes A_{q+1} \otimes \cdots \otimes A_n
\end{equation*}
for $ 1 \le p < q \le n $.  The cactus group is a group with these generators and some natural relations (see section \ref{se:cacgr}).

\subsection{Moduli space of curves}

The cactus group surjects onto $ S_n $ and its kernel is the fundamental group of the Deligne-Mumford compactification, $ \overline{M}_0^{n+1}=\overline{M}_0^{n+1}(\R) $, of the moduli space of real genus 0 curves with $n+1 $ marked points.  The generator $ s_{pq} $ corresponds to a path in $ \overline{M}_0^{n+1} $ in which the marked points $ p, \dots, q $ bubble off onto a new component and then return to the main component in the reversed order.  The cactus group in relation to this fundamental group first appeared in the work of Devadoss \cite{Dev} and Davis-Januszkiewicz-Scott \cite{DJS}.

Elements of this moduli space look somewhat like cacti from the genus Opuntia.  This justifies the name cactus group.  In a similar vein, we would like to propose the name cactus categories for coboundary categories, just as braided categories have replaced the earlier name of quasi-triangular categories.

\subsection{Acknowledgements}
We would like to thank A. Berenstein for explaining to us the construction of the commutor for crystals.  We also benefited from helpful conversations with L. Bartholdi, P. Etingof,  A. Knutson, N. Reshetikhin, N. Snyder, and D. Thurston.  The second author wishes to thank Y. Lai for her hospitality during a visit to Davis and for drawing our attention to the work of S. Devadoss.

\section{Crystals}

A crystal should be thought of as a combinatorial model for a representation of a 
Lie algebra $ \mathfrak{g} $.  

Let $ \mathfrak{g} $ be a complex reductive Lie algebra.  Let $ \Lambda $ denote its 
weight lattice, $\Lambda_+$ denote its set of dominant weights, $ I $ denote the set of vertices of its Dynkin diagram, 
$ \{ \alpha_i \}_{i \in I} $ denote its simple roots, and $ \{ \alpha_i^\vee \}_{i \in I} $ denote
its  simple coroots.
 
We follow the conventions in Joseph \cite{Joseph} in defining crystals, except that we will only consider what he calls ``normal crystals''.  

A {\bf $\mathfrak{g}$-crystal} is a finite set $ B $ along with maps:
\begin{gather*}
\wt : B \rightarrow \Lambda \\
\varepsilon_i, \phi_i : B \rightarrow \mathbb{Z} \\
e_i, f_i : B \rightarrow B \sqcup \{0\} 
\end{gather*}
for each $ i \in I $ such that:
\begin{enumerate}
\item
for all $ b \in B $ we have $ \phi_i(b) - \varepsilon_i(b) = \langle 
\wt(b),\alpha_i^\vee \rangle $ 
\item
$ \varepsilon_i(b) = \max \{ n : e_i^n \cdot b \ne 0 \} $ and $ \phi_i (b) = \max \{ n 
: f_i^n \cdot b \ne 0 \} $ for all $ b \in B $ and $ i \in I $
\item
if $ b \in B $ and $ e_i \cdot b \ne 0 $ then $ \wt(e_i \cdot b) = \wt(b) + 
\alpha_i $, similarly if $ f_i \cdot b \ne 0 $ then $ \wt(f_i \cdot b) = \wt(b) 
- \alpha_i $
\item
for all $ b, b' \in B $ we have $ b' = e_i \cdot b $ if and only if $ b = f_i 
\cdot b' $
\end{enumerate}

We think of $ B $ as the basis for some representation of $ \mathfrak{g} $, with the $e_i $ 
and $ f_i $ representing the actions of the Chevalley generators of $ \mathfrak{g} $.

Let $ A, B $ be crystals.  They have a \textbf{tensor product} $ A \otimes B $ 
defined as follows.  The underlying set is $ A \times B $ and $\wt(a,b) = \wt(a) + \wt(b) $.  We define $ e_i, f_i $ by the following formula:

$$e_i\cdot (a , b)=
\begin{cases}
(e_i \cdot a , b), \quad \text{if}\quad  \varepsilon_i(a) > \phi_i(b)\\
(a , e_i \cdot b),\quad \text{otherwise}
\end{cases}$$
$$f_i\cdot (a, b)=
\begin{cases}
(f_i \cdot a, b), \quad \text{if}\quad  \varepsilon_i(a) \geq \phi_i(b)\\
(a, f_i \cdot b),\quad \text{otherwise}
\end{cases}$$

The \textbf{direct sum} of two crystals is simply defined to be their disjoint union.

\begin{Example}
Consider the crystals which correspond the dual and the symmetric square of the standard representation of $\mathfrak{sl}_3$.  
Then their tensor product is:

\vspace{.3cm}
\centerline{\epsfig{file=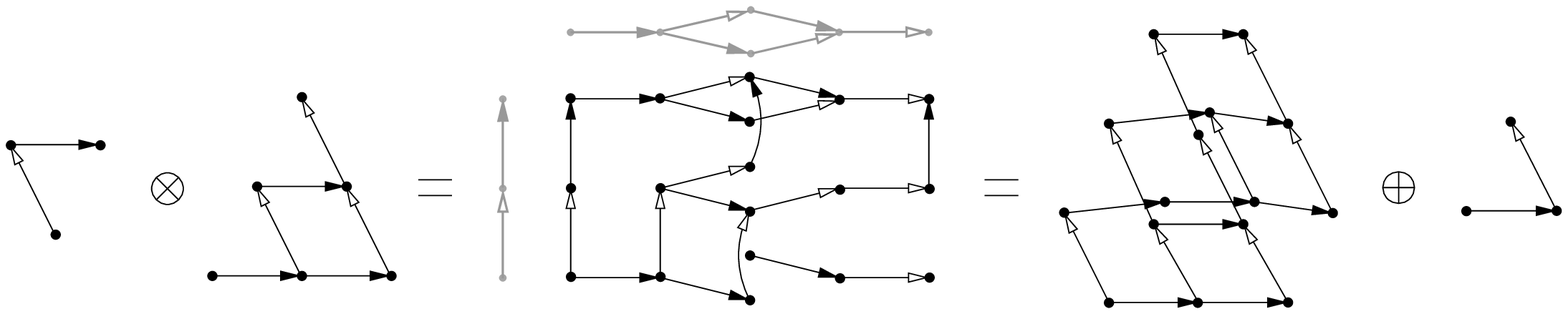,height=2.8cm}}
\vspace{.3cm}
\end{Example}

\subsection{Uniqueness theorem}

We call a crystal {\bf connected} if the underlying graph (where $ b, b' $ are joined by an edge if $ e_i \cdot b = b' $ for some $ i $) is connected.  Similarly we may speak of the {\bf components} of a crystal as the connected components of the underlying graph.  A connected crystal is analogous to an irreducible representation.

A crystal $ B $ is called a {\bf highest weight crystal of highest weight} $ \lambda 
\in \Lambda_+ $, if there exists an element $ b_\lambda \in B $ (called a 
{\bf highest weight element}) that is of weight $ \lambda $, that is sent to $ 0 $ 
by all the $ e_i $, and such that $ B $ is generated by the $ f_i $ acting on $ 
b_\lambda $. 

Not every connected crystal is a highest weight crystal, and there exist non-isomorphic highest weight crystals of the same highest weight.  However, uniqueness does exist for families.

Let $ \mathcal{B} = \{ B_\lambda : \lambda \in \Lambda_+ \} $ be a family of crystals such that $ B_\lambda $ is a highest weight crystal of highest weight $\lambda$.  We say that $ (\mathcal{B}, \iota) $ is a {\bf closed family} if $ \iota_{\lambda, \mu} : B_{\lambda + \mu } \rightarrow B_\lambda \otimes B_\mu $ is an inclusion of crystals.  

\begin{Theorem}[Joseph, \cite{Joseph}]\label{thm:uniqueness}
There exists a unique closed family of crystals $(\mathcal{B}, \iota)$.
\end{Theorem}

There are a number of ways to construct these $ B_\lambda $: 
combinatorially using Littelmann's path model \cite{Lit}, representation theoretically using crystal bases of a quantum group 
representation \cite{kash}, and geometrically using subvarieties of  the affine Grassmannian \cite{BG}.  Each of these construction also produces the inclusions $ \iota $.  If $ \mathfrak{g} = \GLn $, then we can take $ B_\lambda $ to be the set of semi-standard Young Tableaux of shape $ \lambda $.  The weight function is the usual weight of a tableau and the action of $ e_i $ is by the Lascoux-Sch\"utzenberger operator and changes an $ i $ in the tableau to an $ i+1$ \cite{KaNa}.

In addition to being highest weight, these crystals are also lowest weight.
\begin{Proposition} \label{thm:lowestweight}
The crystals $ B_\lambda $ possess a unique lowest weight element $c_\lambda$ of weight $w_0\cdot\lambda$. The $c_\lambda$ are annihilated by all $f_i$ and generate $ B_\lambda $ under the action of the $e_i$.
\end{Proposition}

\begin{proof}
Since the elements of $ B_\lambda $ are partially ordered by their weight and the $ f_i $ decrease the weight, it suffices to prove that there is a unique element annihilated by all the $f_i$.  It is possible to prove this using crystal bases but we will give a purely combinatorial proof.

We appeal to the Lakshmibai-Seshadri chain model.  This is a particular case of the Littelmann path model where we start with a straight line path.  In this model:
\begin{equation*}
B_\lambda = \{ (\mu_0, \dots, \mu_l \in W \cdot \lambda, \, 0 = b_0 < b_1 < \dots < b_l < 1) : \mu_0 <_{b_1} \mu_1 <_{b_2} \dots <_{b_l} \mu_l \}
\end{equation*}
where the $ b_i $ are rational and $ \mu <_b \nu $ denotes the $ b$-Bruhat order.  See Stembridge \cite{stem} for more details.  For our present purposes, all we need is that if $ \mu <_b \nu $, then $ \mu < \nu $.

We identify $ ((\mu_0, \dots, \mu_l), (b_0, \dots, b_l)) $ with the piecewise constant function $ x(t) : (0,1] \rightarrow W \cdot \lambda $ which equals $ \mu_k $ if $ b_{k-1} < t \le b_k $.

By \cite{stem}, we see that:
\begin{equation*}
\phi_i(x(t)) = - \max_{0 \le r \le 1} \int_r^1 \langle x(s), \alpha_i^\vee \rangle 
\end{equation*}

Suppose that $ x(t) $ is killed by all the $ f_i $.  Then $ \phi_i(x(t)) = 0 $ for all $ i$.  Hence for all $ i $, $ \langle \mu_l , \alpha_i^\vee \rangle \le 0 $.  Hence $ \mu_l = w_0 \cdot \lambda $.  Hence $ l = 0 $ (since $\mu_l $ is smaller than any other element in $ W \cdot \lambda $).  So the unique element of $ B_\lambda $ killed by all the $ f_i $ is $ c_\lambda = (w_0 \cdot \lambda, 0) $.  
  
\end{proof}

\subsection{Category of crystals}

The category $\mathfrak{g} $-$ \Crystals$ is 
the category whose objects are crystals $ B $ such that each connected component of $ B $ is isomorphic to some $ B_\lambda $.  For the rest of this paper, crystal means an object in this category.  A 
morphism of crystals is a map of the 
underlying sets that commutes with all the structure maps, including the two 0 elements.  (We might more accurately call our category the category of crystal bases of 
the associated quantum group, since the crystals that arise from crystal bases are exactly those of this form).

If $ \psi : A \rightarrow B $ is a morphism of crystals and $ a \in A $ is killed by all the $e_i$, then $ \psi(a) $ is also killed by all the $ e_i $.  Thus if $ A $ and $ B $ are highest weight crystals, then $ \psi $ must take the highest weight element of $ A $ to the highest weight element of $ B $.  This immediately implies the following version of Schur's Lemma.

\begin{Lemma} 
$ \Hom(B_\lambda, B_\mu) $ contains just the identity if $ \lambda = \mu $ and 
is empty otherwise.  Hence if $ B $ is a crystal there is exactly one way to identify each of its components with a $ B_\lambda $.
\end{Lemma}

The more interesting structure in $ \mathfrak{g}$-$\Crystals$ is the tensor product.  First note that 
the tensor product of crystals is inherently associative: i.e. if $ A, B, C $ are crystals then
\begin{align*}
\alpha_{A,B,C} : (A \otimes B) \otimes C &\rightarrow A \otimes (B \otimes C) \\
((a,b),c) &\mapsto (a,(b,c)) 
\end{align*}
 is an isomorphism.  So we can drop parenthesization when dealing with repeated tensor products.

Note that the definition of the tensor product of crystals is not symmetric, i.e. the map 
\begin{align*} \flip : A \otimes B &\rightarrow  B \otimes A \\  (a,b) &\mapsto (b, a) 
\end{align*}
 is not a morphism of crystals.  

However we will now define isomorphisms $ \sigma_{A, B} : A \otimes B \rightarrow B \otimes A $.  We call such a family a {\bf commutor}.  Following an idea of A. Berenstein (which has not yet appeared in the literature), we produce these isomorphisms by first defining a ``reversing'' involutions $ \xi_B : B \rightarrow B $ for each crystal $ B $.  These will not be a morphisms of crystals, just automorphisms of the underlying graphs which exchange highest weight and lowest weight elements. Such maps were also considered in \cite{Sc}.

Let  $ \theta : I \rightarrow I $ be the Dynkin diagram automorphism such that $\alpha_{\theta(i)} = - w_0 \cdot \alpha_i $ where $ w_0 $ denotes the long element in the Weyl group of $ \mathfrak{g} $.  For $ \mathfrak{g} = \GLn $ with the usual numbering of the simple roots, $ \theta(i) = n- i $.

Let $ \overline{B}_\lambda $ denote the crystal with underlying set $ \{ \overline{b} : b \in B_\lambda \} $ and crystal structure:
\begin{gather*}
e_i \cdot \overline{b} = \overline{f_{\theta(i)} \cdot b} \quad
f_i \cdot \overline{b} = \overline{e_{\theta(i)} \cdot b} \quad
\wt (\overline{b}) =  w_0 \cdot \wt(b)
\end{gather*} 

\begin{Lemma}
With this definition, $ \overline{B}_\lambda $ is a highest weight crystal of highest weight $ \lambda $.  Moreover $ \flip\, \circ\, \iota_{\mu, \lambda} : \overline{B}_{\lambda+\mu}  \rightarrow \overline{B}_\lambda \otimes \overline{B}_\mu $ is an inclusion of crystals, where $ \iota $ is the inclusion map for the original family $ B_\lambda $.
\end{Lemma}

\begin{proof}
The only non-trivial thing to check is that $\overline B_\lambda$ is highest weight. This is the content of Proposition \ref{thm:lowestweight}.

%
%
\end{proof}

By Theorem \ref{thm:uniqueness} and Schur's Lemma, the above lemma gives us a crystal isomorphism $ \overline{B}_\lambda \rightarrow B_\lambda $.  Composing with the map of sets $ B_\lambda \rightarrow \overline{B}_\lambda $ given by $ b \mapsto \overline{b} $, we get a map of sets $\xi= \xi_{B_\lambda} : B_\lambda \rightarrow B_\lambda $.  From the crystal structure on $ \overline{B}_\lambda $ we see that:

\begin{equation}
\begin{gathered} \label{twistrel}
e_i \cdot \xi(b) = \xi(f_{\theta(i)} \cdot b) \\
f_i \cdot \xi(b) = \xi(e_{\theta(i)} \cdot b) \\
\wt (\xi(b)) =  w_0 \cdot \wt(b)
\end{gathered}  
\end{equation}
for all $ b \in B_\lambda $.

In the case  $ \mathfrak{g} = \GLn $, with the tableaux model for $ B_\lambda $,  $ \xi_{B_\lambda} $ is the Sch\"utzenberger involution on tableaux (see \cite{llt}).

We now define $ \xi_B : B \rightarrow B $ for any crystal $ B $ by applying the appropriate $ \xi_{B_\lambda} $ to each connected component.  This map $ \xi $ satisfies the properties in equation (\ref{twistrel}). 

Note that $ \xi_{B_\lambda} \circ \xi_{B_\lambda} $ is a map of crystals, hence by Schur's Lemma $ \xi_{B_\lambda} \circ \xi_{B_\lambda} = 1 $, and so $ \xi_B \circ \xi_B = 1$ for any crystal $B$.

Let $ A, B $ by crystals.  We define:
\begin{align*}
\sigma_{A, B} :  A \otimes B &\rightarrow B \otimes A  \\
(a,b) &\mapsto \xi_{B\otimes A} (\xi_B(b), \xi_A(a))
\end{align*}

\begin{Theorem}
The map $\sigma_{A, B} $ is an isomorphism of crystals and is natural in $A$ and $B$. 
\end{Theorem}

\begin{proof}
Let $ a \in A $ and $ b \in B$.  If $\varepsilon_i(a) > \phi_i(b) $ then $ \varepsilon_{\theta(i)}(\xi(b)) < \phi_{\theta(i)}(\xi(a))$, so :
\begin{equation*}\begin{split}
e_i \cdot \sigma(a,b) &= e_i \cdot \xi(\xi(b), \xi(a)) 
= \xi(f_{\theta(i)} \cdot (\xi(b), \xi(a))) \\
&= \xi(\xi(b), f_{\theta(i)} \cdot \xi(a)) 
= \xi(\xi(b), \xi(e_i \cdot a)) 
= \sigma(e_i \cdot a,b)
= \sigma(e_i \cdot (a,b))
\end{split}\end{equation*}
and similarly for the other case.  So $ \sigma $ commutes with $ e_i $.  Similarly, $ \sigma $ commutes with $ f_i $.  Hence $ \sigma $ is a map of crystals.  

The map $ \sigma $ is a natural isomorphism since both $\xi$ and $ \flip $ are bijections which are natural with respect to maps of crystals (more precisely, $\xi$ is an automorphism of the forgetful functor forget:$\mathfrak{g} $-$\Crystals \to \mathsf{Sets}$, and $\flip$ is an automorphism of forget$\circ \otimes$).  
\end{proof}

Note that $ \sigma_{A,B} = \xi_{B \otimes A} \circ (\xi_B \otimes \xi_A) \circ \flip$ .  This suggests that we consider $ \flip \circ\, (\xi_A \otimes \xi_B) \circ \xi_{A \otimes B} $.  In fact we have:

\begin{Proposition}\label{pro:2sigmas}
$ \sigma_{A,B} =  \flip \circ\, (\xi_A \otimes \xi_B) \circ \xi_{A \otimes B} $.  
\end{Proposition}

\begin{proof}
First we note that since $ \flip \circ\, (\xi_A \otimes \xi_B) \circ \xi_{A \otimes B}=\sigma_{B,A}^{-1} $, it is also an isomorphism of crystals.  

Let $ (a,b) \in A \otimes B $.  Let $(c,d) = \xi(a,b) $.  Then we want to show that $ (\xi(d), \xi(c)) = \xi(\xi(b), \xi(a)) $.  Since both maps are isomorphisms of crystals, by Schur's Lemma it suffices to show that $ (\xi(d), \xi(c)) $ lies in the same component as $ \xi(\xi(b), \xi(a)) $.  Now by construction $ \xi $ preserves components so it suffices to check that $ (\xi(b), \xi(a)) $ lies in the same component as $ (\xi(d), \xi(c)) $.  But $ (c,d) $ lies in the same component as $ (a,b) $, hence there exist some $ i_1, \dots i_r, j_1, \dots j_s $ such that $ e_{i_1}  \cdots e_{i_r} \cdot f_{j_1} \cdots f_{j_s} \cdot (a,b) = (c, d) $.  But then it is easily checked that $ f_{\theta(i_1)}  \cdots  f_{\theta(i_r)} \cdot e_{\theta(j_1)}  \cdots e_{\theta(j_s)} \cdot (\xi(b), \xi(a)) = (\xi(d), \xi(c)) $.  Hence the two elements lie in the same component and so the result follows.
\end{proof}

\begin{Example}
For the tensor product considered in Example 1, the map $\sigma $ looks like

\vspace{.3cm}
\centerline{\epsfig{file=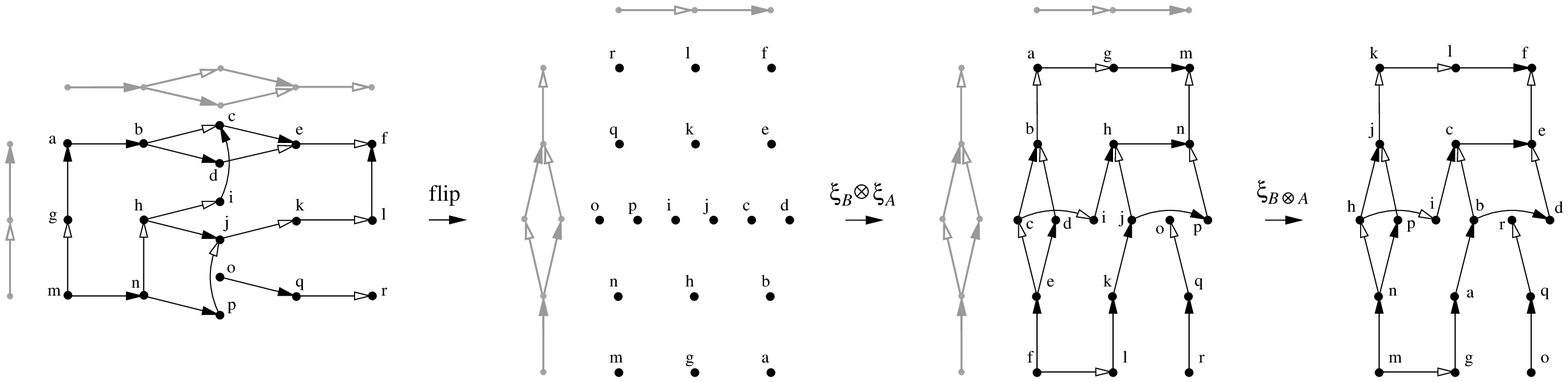,height=3.6cm}}
\vspace{.3cm}
\end{Example}

\subsection{Properties of the commutor} 
Now that we have a canonically defined commutor for the category of crystals, it is interesting to consider what axioms it obeys. First we have symmetry:

\begin{Proposition}
\begin{equation*}
 \sigma_{B,A} \circ \sigma_{A, B} = 1
\end{equation*}
\end{Proposition}
\begin{proof}
This follows directly from $ \xi \circ \xi = 1 $ and Proposition \ref{pro:2sigmas}.
\end{proof}

Now we consider the compatibility between the commutor and the associator.   The commutor does not obey the usual hexagon\footnote{If we had included the associator isomorphisms, then (\ref{eq:hex}) would have the shape of an hexagon.} axiom:
\begin{equation} \label{eq:hex}
\xymatrix{
A \otimes B \otimes C \ar[rr]^{\sigma_{A \otimes B, C}} \ar[dr]_{1 \otimes \sigma_{B, C}} & & C \otimes A \otimes B \\
& A \otimes C \otimes B \ar[ur]_{\sigma_{A, C} \otimes 1} &  \\
}
\end{equation}

In fact even for $ \mathfrak{sl}_2$, it is not possible to find natural isomorphisms $ \beta_{A,B} : A \otimes B \rightarrow B \otimes A $ which obey the hexagon axiom.

Suppose $ \beta $ is such a family of isomorphisms.  Let $ A = \{a_0, a_1 \}, B = \{b_0, b_1 \}, C= \{c_0, c_1 \} $ denote three copies of the crystal of the standard representation of $ \mathfrak{sl}_2 $.  Let $ D = \{d\} $ denote the crystal of the trivial representation.  

We have an inclusion $ j : D \rightarrow A \otimes B $ which takes $ d $ to $ (a_0, b_1) $.  Since $ \beta $ is natural, the diagram
\begin{equation*}
\xymatrix{
D \otimes C \ar[d]_{j \otimes 1} \ar[r]^{\beta_{D, C}} & C \otimes D \ar[d]^{1 \otimes j}\\
(A \otimes B) \otimes C \ar[r]_{\beta_{A \otimes B, C}} & C \otimes (A \otimes B) \\
}
\end{equation*}
commutes.  This shows that $ \beta_{A \otimes B, C}(a_0, b_1, c_0) = (c_0, a_0, b_1) $.

Furthermore, $ \beta_{B,C} (b_1, c_0) = (c_1, b_0) $, because there is only one isomorphism between $ B \otimes C $ and $ C \otimes B $.  
Similarly, $ \beta_{A, B}(a_0, c_1) = (c_0, a_1) $ .  Thus, $ (\beta \otimes 1) \circ (1 \otimes \beta) (a_0, b_1, c_0) = (c_0, a_1, b_0) $.

Hence the two arcs of the ``hexagon'' (\ref{eq:hex}) give different results when applied to $ (a_0, b_1, c_0) $ and so we see that the hexagon does not commute.

Returning to our commuter $ \sigma $, we do have the following compatibility result.

\begin{Theorem}
The following diagram commutes in $ \mathfrak{g}$-$Crystals$:
\begin{equation} \label{eq:bubc}
\xymatrix{
A \otimes B \otimes C \ar[d]_{\sigma_{A,B} \otimes 1} \ar[r]^{1 \otimes \sigma_{ B, C}} & A \otimes C \otimes B \ar[d]^{\sigma_{A, C \otimes B}}\\
B \otimes A \otimes C \ar[r]_{\sigma_{B \otimes A, C}} & C \otimes B \otimes A \\
}
\end{equation}
\end{Theorem}

Note that by naturality of $ \sigma $, we have $ \sigma_{B \otimes A, C} \circ ( \sigma_{A,B} \otimes 1) = (1 \otimes \sigma_{A,B}) \circ \sigma_{A \otimes B, C}$ , so the above diagrams implies three other diagrams of the same form.

\begin{proof}
Let $ (a, b, c) \in A \otimes B \otimes C $.\\
Following down and right gives:
\begin{equation*}
 (a,b,c) \mapsto (\xi(\xi(b), \xi(a)), c) \mapsto \xi(\xi(c), \xi(b), \xi(a))
\end{equation*}
Following right and down gives:
\begin{equation*}
(a,b,c) \mapsto (a, \xi(\xi(c), \xi(b))) \mapsto \xi(\xi(c), \xi(b), \xi(a))
\end{equation*}
\end{proof}

\subsection{A commutor for quantum groups}
It is natural at this point to ask how the commutor defined above is related to the associated quantum group. 

Let $ U_q(\mathfrak{g}) $ denote the quantized universal enveloping algebra of the Lie algebra of $ \mathfrak{g}$.  Recall that this is a Hopf algebra over $ \mathbb{C}(q) $ with generators $ E_i, F_i $ for $ i \in I $ and $ K_\alpha $ for $ \alpha$ in the coroot lattice and the following relations:
\begin{gather*}
K_\alpha K_\beta = K_{\alpha + \beta} \\
K_\alpha E_i K_{\alpha}^{-1} = q_i^{\langle \alpha_i, \alpha \rangle} E_i \\
K_\alpha F_i K_{\alpha}^{-1} = q_i^{-\langle \alpha_i, \alpha \rangle} F_i \\
E_i F_i - F_i E_i = \frac{K_i - K_i^{-1}}{q_i - q_i^{-1}}
\end{gather*}
where $ q_i = q^{d_i}$ (where $d_i$ are coprime integers chosen to symmetrize the Cartan matrix).  We also have the quantum Serre relations which we will not pause to write down.  As is standard, we write $ K_i $ for $ K_{\alpha_i^\vee}$.  For more information, see \cite{CP}.

The coproduct structure on $ U_q(\mathfrak{g}) $ is given by:
\begin{gather*}
\Delta(K_\alpha) = K_\alpha \otimes K_\alpha \\
\Delta(E_i) = E_i \otimes K_i + 1 \otimes E_i \\
\Delta(F_i) = F_i \otimes 1 + K_i^{-1} \otimes F_i^{-1} 
\end{gather*}

We will construct a commutor for the category of $\qg $ modules by a procedure analogous to our construction for crystals.

Define an algebra automorphism $ \xi : \qg \rightarrow \qg $ by:
\begin{equation*}
\xi(E_i) = F_{\theta(i)}, \quad \xi(F_i) = E_{\theta(i)}, \quad \xi(K_\alpha) = K_{w_0 \cdot \alpha}
\end{equation*}

By examining the relations for the product and the definition of the coproduct above, we see that:

\begin{Proposition} 
These formulas extend to an algebra automorphism with $ \xi \circ \xi = 1 $.  Moreover, $ \xi $ is a coalgebra antiautomorphism, i.e.:
\begin{equation*}
(\xi \otimes \xi)( \Delta^{op}(a)) = \Delta(\xi(a))
\end{equation*}
\end{Proposition}

Now, suppose that $ V $ is a finite-dimensional $ \qg $ module.  Define a new $ \qg $ module $\overline{V} $ to be the same underlying vector space, written as $ \{ \overline{v} : v \in V \} $, but with the action twisted by $ \xi $:
\begin{equation*}
a \cdot \overline{v} = \overline{\xi(a) \cdot v}
\end{equation*}

Since $ \xi $ acts on the $ K_\alpha $ by $ w_0 $, we see that the character of $ \overline{V} $ is $ w_0 $ applied to the character of $ V $.  But the character of $ V $ is invariant under $ w_0 $, so $ \overline{V} $ and $ V $ have the same character.  Hence they are isomorphic.  In particular, if $ V = V_\lambda $ is irreducible, $ V_\lambda $ and $ \overline{V_\lambda} $ are isomorphic by an isomorphism that is unique up to scalar.  

Hence we can fix an automorphism of vector spaces $ \xi_{V_\lambda} : V_\lambda \rightarrow V_\lambda $ such that:
\begin{equation} \label{xiprop}
\xi(a) \cdot \xi_{V_\lambda} (v) = \xi_{V_\lambda} (a \cdot v) \ \text{ for } a \in \qg \text{ and } v \in V_\lambda
\end{equation}

Also, we may choose the normalization of $ \xi_{V_\lambda} $ such that $ \xi_{V_\lambda} \circ \xi_{V_\lambda} = 1 $.  This makes $ \xi_{V_\lambda} $ unique up to sign. We will show in Theorem \ref{prescan} that there exists a preferred sign convention for $\xi_{V_{\lambda}}$ that makes it compatible with the canonical basis.

These $ \xi_{V_\lambda} $ patch together to form a map $ \xi_V : V \rightarrow V  $ for any representation $ V $, defined by $ \xi_V(\phi(w)) = \phi(\xi_{V_\lambda}(w)) $ where $ \phi : V_\lambda \rightarrow V $ is any morphism of $ \qg $ modules.  The resulting $ \xi_V $ also satisfies (\ref{xiprop}) and $ \xi_V \circ \xi_V = 1 $.

Now we define a commutor in the same way as for the crystals. If $ V$ and $ W $ are $ \qg $ modules, we let:
\begin{equation}\begin{split}\label{sigmaVW}
\sigma_{V, W} : V \otimes W & \rightarrow W \otimes V \\
v \otimes w &\mapsto \xi_{W\otimes V}(\xi(w) \otimes \xi(v))
\end{split}\end{equation}

We have:

\begin{Theorem} \label{quansig}
\begin{enumerate}
\item $ \sigma_{V,W} $ is an isomorphism of $\qg$ modules and it is natural in $V$ and $W$.
\item It obeys the cactus axiom (\ref{eq:bubc}) 
\item $ \sigma_{V,W} =  \flip \circ\, (\xi_V \otimes \xi_W)  \circ \xi_{V,W} $
\item $ \sigma_{W,V} \circ \sigma_{V,W} = 1 $
\end{enumerate}
\end{Theorem}

\begin{proof}
\begin{enumerate}
\item
Let $ a \in \qg $ and $ v \in V \otimes W $.  Then:
\begin{align*}
\xi(\xi \otimes \xi (\flip(a \cdot v))) &= \xi( \xi \otimes \xi( \flip(\Delta(a) \cdot v))) \\
&= \xi( \xi \otimes \xi( \Delta^{op}(a) \cdot \flip (v))) \\
&= \xi( \xi \otimes \xi(\Delta^{op}(a)) \cdot \xi \otimes \xi (\flip(v))) \\
&= \xi( \Delta(\xi(a)) \cdot \xi \otimes \xi (\flip(v))) \\
&= \xi( \xi(a) \cdot \xi \otimes \xi (\flip(v))) \\
&= a \cdot \xi( \xi \otimes \xi (\flip(v))) \\
\end{align*}

\item
This proof is identical to the proof for the commutor in $\mathsf{Crystals}$.

\item
Let $ v \in V \otimes W $ be in the image of some $ \phi : V_\lambda \rightarrow V \otimes W $. Since any two elements of $ V_\lambda $ are related by the action of $\qg$, we can write $ \xi(v) = Y \cdot v $ for some $ Y \in \qg $.  Then $ \xi(Y) \cdot \xi(v) = v $. Since $ \sigma $ is a $\qg$ module morphism and $\xi$ commutes with $\qg$ module morphisms, we also have $ \xi(Y) \cdot \xi(\sigma(v)) = \sigma(v) $ .
So:
\begin{align*}
 \flip \circ \,(\xi \otimes \xi) \circ \xi (v) &= \xi \otimes \xi ( \flip(Y \cdot v)) \\ 
&= \xi \otimes \xi (\Delta^{op}(Y)\cdot (\flip(v))) \\
&= \xi(Y) \cdot \xi \otimes \xi (\flip(v)) \\
&= \xi(Y) \cdot \xi(\sigma(v)) = \sigma(v)
\end{align*}

\item
The symmetry $ \sigma_{V, W} \circ \sigma_{W, V} = 1 $ follows easily form (iii) since $ \xi^{-1} = \xi $.

\end{enumerate}
\end{proof}

\subsection{Relation between commutors}
This commutor is related to the commutor for crystals.  Recall that by the work of Lusztig, each irreducible $ \qg $ module $ V_\lambda $ has a canonical basis $ B $ such that the action of the operators $ E_i, F_i \in \qg $ induces (is a somewhat complicated way) a crystal structure on $ B $ (see \cite{Lus} for more details).

This automorphism $ \xi $ was first considered in the $ \mathfrak{gl}_n$ case by Berenstein-Zelevinsky \cite{BZ}.  They proved the following result in that case which was the inspiration for our more general construction:

\begin{Theorem} \label{prescan}
Let $ V_\lambda $ be the irreducible $\qg$ module of highest weight $\lambda $ and let $ B $ be its canonical basis.  Then there exists a choice of sign normalization for $ \xi_{V_\lambda} $ which preserves $ B $ and agrees with the crystal map $ \xi_{B_\lambda} $.
\end{Theorem}

\begin{proof}
Consider the Hopf algebra automorphism $ \omega: E_i \mapsto F_i, F_i \mapsto E_i, K_i \mapsto K_i^{-1}$.  If we consider the twisted module structure $ V_\lambda^{\omega}$ induced by $\omega $, we see that $ V_\lambda $ is isomorphic to $ V_{\theta( \lambda)} $ where $ \theta(\lambda) = -w_0 \cdot \lambda $.  In \cite[21.1.2]{Lus}, Lusztig shows that there exists a choice of this isomorphism $\omega_\lambda : V_\lambda \to V_{\theta(\lambda)}$ that sends the canonical basis of $ V_\lambda $ to the canonical basis of $ V_{\theta(\lambda)}$.  

Our automorphism $ \xi $ is the composition of $\omega$ with the automorphism of $ \qg $ induced from the Dynkin diagram automorphism $ \theta $.   This automorphism $ \theta $ induces an isomorphism of representations $ \theta_\lambda : V_\lambda  \rightarrow V_{\theta (\lambda)} $ which takes the canonical basis to the canonical basis.

Hence, we can choose $ \theta\circ\omega:V_\lambda\to V_{\theta(\lambda)}\to V_\lambda$ to be our $\xi_{V_\lambda}$ since this composition satisfies the defining property (\ref{xiprop}).

 And so $ \xi_{V_\lambda} : V_\lambda \rightarrow V_\lambda $ will preserve the canonical basis.  Since $ \xi $ exchanges $ E_i $ and $ F_{\theta(i)} $, $ \xi_{V_\lambda} $ will give a map on $ B $ that satisfies the same properties (\ref{twistrel}) as the crystal map $ \xi_{B_\lambda} $.  Since $ \xi_{B_\lambda} $ is the unique map satisfying these properties, the result follows.
\end{proof}

\subsection{Relation with braiding} \label{se:drin}
Recall that there is a more standard commutor (usually called the braiding) on the category of $ \qg $ representations which is constructed using the universal $ R$-matrix for $ \qg$ (see \cite{CP} for more details).  This element $ R $ lives in $ U_h(\mathfrak{g}) \otimes U_h(\mathfrak{g}) $ where $ U_h(\mathfrak{g}) $ is the formal form of the quantum group.  The braiding is given by:
\begin{align*}
B : v \otimes w \mapsto R \cdot w \otimes v
\end{align*}

In \cite{Drin}, Drinfel'd explained how use $ R $ to construct another commutor which is symmetric and satisfies the cactus axiom (\ref{eq:bubc}).  Namely set $ R' = R(R^{21}R)^{-1/2} $ where the square root is taken with respect to the ``$h$'' filtration on $  U_h(\mathfrak{g}) \otimes U_h(\mathfrak{g}) $ and where $ R^{21} = \flip(R) $ .  Then the new commutor is defined by:
\begin{align*}  
\sigma' : v \otimes w \mapsto R' \cdot w \otimes v
\end{align*}
Drinfel'd calls this process \textbf{unitarization}.

The following question was suggested by A. Berenstein:
Does there exist an appropriate choice of normalization for which our commutor $ \sigma $ agrees with Drinfel'd's $ \sigma'$?

The normalization is quite subtle as we have a choice of $ \pm 1$ normalization for each isotypic component of each $ V_\lambda \otimes V_\mu $.

\section{Coboundary categories}
We now put the category of crystals in a more general framework and investigate its structure.

A {\bf monoidal category} is a category $ \Cat $ along with a functor $ \otimes : \Cat \times \Cat \rightarrow \Cat $ and natural isomorphisms (called the {\bf associator}) $ \alpha_{A, B, C} : A \otimes (B \otimes C) \rightarrow (A \otimes B) \otimes C $ such that the following pentagon diagram commutes:
\begin{equation*}
\xymatrix{
A \otimes (B \otimes (C \otimes D)) \ar[r] \ar[d] & (A \otimes B) \otimes (C \otimes D) \ar[r] & ((A \otimes B) \otimes C) \otimes D \ar[d] \\
A \otimes ((B \otimes C) \otimes D) \ar[rr] & & (A \otimes (B \otimes C)) \otimes D) \\
}
\end{equation*}

There is also a unit object $1\in\Cat$ and isomorphisms $\eta_R:1\otimes A\to A \leftarrow A\otimes 1:\eta_L$ such that the composites $A\to 1\otimes A\to A$ and $A\to A\otimes 1 \to A$ are the identity. The unit will not play an important role in our discussion.

A {\bf coboundary category} \cite{Drin} is a monoidal category along with natural isomorphisms $ \sigma_{A, B} : A \otimes B \rightarrow B \otimes A $ (called the {\bf commutor}) such that $ \sigma_{B,A} \circ \sigma_{A,B} = 1 $ and the following diagram commutes:
\begin{equation*}
\xymatrix{
A \otimes (B \otimes C) \ar[r]^{1\otimes\sigma_{B,C}} \ar[d]_{\alpha_{A, B, C}} &
A \otimes (C \otimes B) \ar[r]^{\sigma_{A, C \otimes B}} &
(C \otimes B) \otimes A \\
(A \otimes B) \otimes C \ar[r]^{\sigma_{A,B}\otimes 1} &
(B\otimes A)\otimes C \ar[r]^{\sigma_{B \otimes A,C}} &
C \otimes (B \otimes A) \ar[u]_{\alpha_{C,B,A}} \\
}
\end{equation*}
We call these two conditions the symmetry axiom and the cactus axiom. We also ask that $\eta_L^{-1} \circ \eta_R=\sigma_{1,A}$.

This should be compared with the definition of a braided monoidal category, where the commutor is not required to be symmetric and the compatibility between the associator and the commutor is expressed by two hexagon diagrams that involve three commutor moves and three associator moves (for a picture without the associators, see (\ref{eq:hex})).

Our results of the previous section can now be restated as:
\begin{Theorem}
$\mathfrak{g}$-$\Crystals$ and $ \qg$-$\mathsf{Modules} $ form coboundary categories using the above commutors.
\end{Theorem}

\subsection{Cactus group action} \label{se:cacgr}

We now want to consider all morphisms of iterated tensor products $ A_1 \otimes \dots \otimes A_n $ that can be made using our commutor.  For example, $ 1 \otimes \sigma_{B \otimes C, D \otimes E } : A \otimes B \otimes C \otimes D \otimes E \rightarrow A \otimes D \otimes E \otimes B \otimes C $ (we can drop parenthesization because the pentagon axiom ensures that there is a unique isomorphism between any two parenthesizations).

For braided categories, the hexagon axiom ensures that all such morphisms can be written as compositions of switches of adjacent factors.  It also shows that switching adjacent factors obeys the Yang-Baxter equations:
\begin{equation*}
\sigma_{12} \circ \sigma_{23} \circ \sigma_{12} = \sigma_{23} \circ \sigma_{12}\circ \sigma_{23}
\end{equation*}
where $ \sigma_{12} $ denotes the morphism:
$ A \otimes B \otimes C \rightarrow B \otimes A \otimes C $.

These results imply that the braid group acts on multiple tensor products and that this action contains all the morphisms generated by the commutor.  More precisely, let $ \Cat $ be a braided category and let $ A_1 \dots A_n \in \Cat$.  If $ \rho \in B_n $ (the braid group on $n $ strands), then we get a morphism in $ \Cat $:
\begin{equation*}
A_1 \otimes \dots \otimes A_n \rightarrow A_{\rho(1)} \otimes \dots \otimes A_{\rho(n)}
\end{equation*}
by expressing $ \rho $ as a product of generators of the braid group and then using the generator to switch adjacent pairs in the tensor product (here $ \rho(i) $ is the standard action of $ B_n $ on $ \{ 1 \dots n \} $ using the usual map $ B_n \rightarrow S_n $).

For coboundary categories, the situation is different but analogous.
To save space we will drop the symbol $ \otimes $ and use the convention that tensor product of objects is denoted by concatenation.

Let $ \Cat $ be a coboundary category and let $ A_1, \dots, A_n \in \Cat $.  
If $ 1 \le p \le r < q \le n $, we get a natural isomorphisms denoted $ \sigma_{p,r,q} $ defined by:
\begin{gather*}
 (\sigma_{p,r,q})_{A_1, \dots, A_n} := 1 \otimes \sigma_{A_p \cdots A_r, A_{r+1} \cdots A_q} \otimes 1:\\ 
A_1 \cdots A_{p-1}A_p \cdots A_r A_{r+1} \cdots A_q A_{q+1}\cdots  A_n \longrightarrow
A_1 \cdots A_{p-1}A_{r+1} \cdots A_q A_p \cdots A_rA_{q+1} \cdots A_n.
\end{gather*}

We are interested in the natural isomorphims generated using these $ \sigma_{p,r,q} $.  The basic building blocks in coboundary categories will be maps
$$s_{p,q}:A_1  \cdots A_{p-1}A_p A_{p+1} \cdots A_{q-1}A_q A_{q+1} \cdots  A_n \longrightarrow A_1  \cdots A_{p-1}A_q A_{q-1} \cdots A_{p+1}A_p A_{q+1} \cdots  A_n$$ which `reverse intervals'. They are defined recursively by $ s_{p, p+1}= \sigma_{p,p,p+1} $, and $ s_{p,q}= \sigma_{p,p,q} \circ s_{p+1, q}$ for $ q-p > 1 $:
\begin{equation*}
A_1  \cdots A_p \cdots  A_q \cdots  A_n \xrightarrow{s_{p+1,q}}
A_1 \cdots A_p A_q \cdots A_{p+1} \cdots  A_n \xrightarrow{\sigma_{p,p,q}}
A_1  \cdots A_q \cdots A_p \cdots  A_n.
\end{equation*}

\noindent By convention, we let $s_{p,p}=1$.
For the commutors constructed in Section 2, we have $$s_{p,q}:(a_1,\ldots, a_n)\mapsto(a_1,\ldots ,a_{p-1},\xi(\xi(a_q)\ldots\xi(a_p)),a_{q+1},\ldots, a_n)$$ and the statements of Lemmas \ref{alterdef} and \ref{lem:rels} can be checked directly. Since we want to show them in an arbitrary coboundary category, we need to work a little bit harder.

\begin{Lemma} \label{alterdef}

\begin{enumerate}
\item $s_{p,q} \circ s_{k,l}=s_{k,l}\circ s_{p,q}\quad$ if $p<q<k<l$.
\item $\sigma_{p,q,r}\circ s_{k,l}=s_{k+q-r,l+q-r}\circ\sigma_{p,q,r}\quad$ if $p\le k < l \le r < q$, \hfill\break
      $\sigma_{p,q,r}\circ s_{k,l}=s_{k-(r+1-p),l-(r+1-p)}\circ\sigma_{p,q,r}\quad$ if $p < r+1 \le k < l \le q$.
\item $s_{p,q}=\sigma_{p, r, q} \circ s_{p,r} \circ s_{r+1, q}$.

\item $ s_{p,q} \circ s_{p,q} = 1 $.

\item Every $ \sigma_{p,r,q} $ can be written as a composition of the $ s_{k,l}$.
\end{enumerate}
\end{Lemma}

\begin{proof}
\begin{enumerate}

\item In a monoidal category $f\otimes 1$ commutes with $1\otimes g$.

\item We expand $s_{k,l}$ in terms of $\sigma$ using the definition, then use naturality to pass the $\sigma$'s across $\sigma_{p,q,r}$. Each time something passes across, we add $q-r$ (or subtract $r+1-p$) to all the indices.

\item 
We use a double induction on $p-q $ and on $ r$.  The base case $r = p $ is the definition.  Assume that $ r > p $. Then:
\begin{equation*}\begin{split}
\sigma_{p,r,q} \circ s_{p,r} \circ s_{r+1,q} &= \sigma_{p,r,q} \circ \sigma_{p,r-1,r} \circ s_{p,r-1} \circ s_{r+1, q} \\
&= \sigma_{p,r-1,q} \circ \sigma_{r,r,q} \circ s_{r+1, q}  \circ s_{p,r-1} \\
&= \sigma_{p,r-1,q} \circ s_{r,q} \circ s_{p,r-1} = s_{p,q}.
\end{split}\end{equation*}
The first equality uses the induction on $p-q$ and the convention $s_{r,r}=1$. The second equality uses the cactus axiom and (i). The third equality is the definition of $s$ and the last equality uses (i) and induction on $r$.

\item
The $ r = p $ and $ r = q - 1 $ cases of (iii), followed by the second case of (ii) give
\begin{equation*}\begin{split}
s_{p,q} \circ s_{p,q} &= \sigma_{p, p, q} \circ s_{p+1,q} \circ \sigma_{p, q-1, q} \circ s_{p, q-1}  \\
&=s_{p, q-1} \circ \sigma_{p, p, q} \circ \sigma_{p, q-1, q} \circ s_{p, q-1} 
\end{split}\end{equation*}
So the symmetry property $ \sigma_{p,p,q} \circ \sigma_{p, q-1, q} = 1 $ implies (iv) by induction on $ q-p $.

\item
Follows from (iii) and (iv).
\end{enumerate}
\end{proof}

We now investigate the relations satisfied by these $ s_{p,q} $.  For $ p < q $ let $ \widehat{s}_{p,q} $ denote the involutive element of the symmetric group $S_n $: 
\begin{equation*}
\widehat{s}_{p,q} = \begin{pmatrix}
1 & \cdots & p -1 & p &\cdots& q & q+1 &\cdots& n \\ 
1 & \cdots& p-1 & q& \cdots &p & q+1 &\cdots &n   
\end{pmatrix}
\end{equation*}
So $ s_{p,q} $ is a natural isomorphism from $ F_1 $ to $ F_{\widehat{s}_{p,q}} $, where $F_\pi:\Cat^n\to\Cat$ denotes the functor $(A_1\ldots A_n)\mapsto A_{\pi(1)}\otimes\cdots\otimes A_{\pi(n)}$, for $\pi\in S_n$.

We say that $ p < q $ and $ k < l $ are {\bf disjoint} if $ q < k $ or $ l < p $. We say that $ p < q$ {\bf contains} $  k < l $ if $ p \le k < l \le q $. 

\begin{Lemma}\label{lem:rels}
If $ p < q $ contains $ k < l $, then
\begin{equation}
s_{p,q} \circ s_{k,l} = s_{m, n} \circ s_{p,q} \quad \text{where } m = \widehat{s}_{p,q}(l), n = \widehat{s}_{p,q}(k) 
\end{equation}
\end{Lemma}

\begin{proof}
First compute $m=p+q-l$ and $n=p+q-k$.
Using (i), (ii) and (iii) of Proposition \ref{alterdef}, we have:
\begin{equation*}\begin{split}
s_{p,q} &= \sigma_{p,l,q} \circ s_{p,l}                                         \circ s_{l+1, q} \\
        &= \sigma_{p,l,q} \circ \sigma_{p,k-1,l}  \circ s_{p,k-1} \circ s_{k,l} \circ s_{l+1, q} \\
        &= (\sigma_{p,l,q} \circ \sigma_{p,k-1,l} \circ s_{p,k-1} \circ s_{l+1, q}) \circ s_{k,l}\quad (*) \\
        &= \sigma_{p,l,q} \circ s_{p,l+p-k} \circ \sigma_{p,k-1,l}  \circ s_{p,k-1}  \circ s_{l+1, q} \\
        &= s_{m,n} \circ (\sigma_{p,l,q} \circ \sigma_{p,k-1,l}  \circ s_{p,k-1}  \circ s_{l+1, q}) \quad (*). \\
\end{split}\end{equation*}
Using the two expressions $(*)$ and (iv) of Proposition \ref{alterdef}, we deduce
\begin{equation*}
s_{p,q}\circ s_{k,l}=\sigma_{p,l,q} \circ \sigma_{p,k-1,l}  \circ s_{p,k-1}  \circ s_{l+1, q}= s_{m, n} \circ s_{p,q}.
\end{equation*}
\end{proof}

Let $ J_n $ be the group with generators $ s_{p, q} $ for $1 \le p < q \le n $ and relations:
\begin{enumerate}
\item $ s_{p,q}^2 = 1 $.
\item $ s_{p,q}  s_{k,l} = s_{k,l} s_{p,q}$ if $ p < q $ and $ k < l $ are disjoint.
\item $ s_{p,q}  s_{k,l} = s_{m, n} s_{p,q} $ if $ p< q $ contains $ k < l $, where $ m = \widehat{s}_{p,q}(l)$ and $ n = \widehat{s}_{p,q}(k)$.
\end{enumerate}
This group has appeared in \cite{Dev} under the name quasi-braid group and it is one of the main examples of what \cite{DJS} call mock reflection groups. We like to call $J_n$ the \bf n-fruit cactus group \rm.

The cactus group $ J_n $ admits a natural map to $ S_n $, $ \rho \mapsto \widehat{\rho}$, extending $ s_{p,q} \mapsto \widehat{s}_{p,q} $.
So we have proven:

\begin{Theorem}\label{cactusact}
Let $ \Cat $ be a coboundary category and let $ A_1 \dots A_n \in \Cat $.  
If $ \rho \in J_n $, we have a natural isomorphism $\tau(\rho;A_1,\ldots,A_n): A_1 \cdots A_n \to A_{\widehat{\rho}(1)} \cdots A_{\widehat{\rho}(n)} $ satisfying
$$
\tau(\rho';A_{\widehat{\rho}(1)},\ldots,A_{\widehat{\rho}(n)})\circ\tau(\rho;A_1,\ldots,A_n)=\tau(\rho\rho';A_1,\ldots,A_n).
$$
These natural isomorphisms are exactly those which can be generated using the commutor.
\end{Theorem}

\subsection{Moduli space of curves}
There is a nice geometric interpretation of $ J_n $ and its action on coboundary categories.  
This geometry concerns the total space $\widetilde{M}_0^{n+1}$ of a line bundle over the Deligne-Knudson-Mumford moduli space $\overline{M}_0^{n+1}= \overline{M}_0^{n+1}(\R) $ of stable real curves of genus $ 0 $ with $ n+1 $ marked points.
Our definitions follow Kapranov \cite{Kap}.  

An {\bf n-fruit cactus} is an algebraic curve $ C $ over $ \mathbb{R} $ along with $n+1 $ distinct marked points of C 
and a non-zero tangent vector to $ C $ at the $( n+1 )$st marked point (which will be called the base point of $ C $) such that:
\begin{enumerate}
\item every irreducible component of $ C $ is isomorphic to $ \mathbb{RP}^1$.
\item $C$ is smooth at the marked points.
\item C has only ordinary double points (normal crossings).
\item The graph of components of $ C $ is a tree.
\item The automorphism group of $ C $ is trivial.  
This means that on each component of $ C $ there are at least three points which are either marked or double, 
except for the base point component which is allowed to have only two.
\end{enumerate}

\begin{figure}[htbp]
  \begin{center}
   \epsfig{file=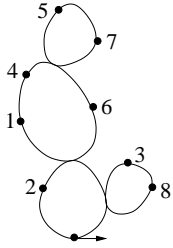,height=3.2cm}
    \caption{An 8 fruited cactus.}
    \label{fig:cactus}
  \end{center}
\end{figure}

One should imagine cacti as being 3-dimensional objects, and the `leaves' being allowed to rotate about their attaching point.

Let $\widetilde{M}_0^{n+1}$ denote the  moduli space of $n$-fruit cacti (the $0$ stands for `genus zero').
It is the total space of the $(n+1)$st tautological line bundle of $\overline M_0^{n+1}$.
The map $p:\widetilde M_0^{n+1}\to \overline M_0^{n+1}$ forgets the vector on the base point and collapses
its component if necessary (i.e. if the base point component has exactly two special points).
The cacti for which this collapsing occurs form the zero section of $p$.

Note that $\widetilde{M}_0^{n+1}$ contains a dense open subset $ U $ consisting of those curves with only one component.  If we map the base point to $ \infty $ and fix the tangent vector, then we see that:
\begin{equation*}
U \cong (\mathbb{R}^n - \Delta)/\mathbb{R}
\end{equation*}
where $ \Delta = \{ (x_1, \dots, x_n) : x_i = x_j \text{ for some } i \ne j \} $ denotes the ``thick diagonal'' and $ \mathbb{R} $ is acting by simultaneous translation.

This open dense subset has $ n! $ different connected components corresponding to the different ways ordering the first $n $ marked points around the curve (or equivalently, the $ n! $ regions of the hyperplane arrangement $ \Delta $).  Note that $ S_n $ acts on $ \widetilde{M}_0^{n+1} $ by permuting the labels of the first $n$ points and that this action is transitive on the $ n! $ connected components of $U\subset \widetilde{M}_0^{n+1}$.

\begin{figure}[htbp]
  \begin{center}
   \epsfig{file=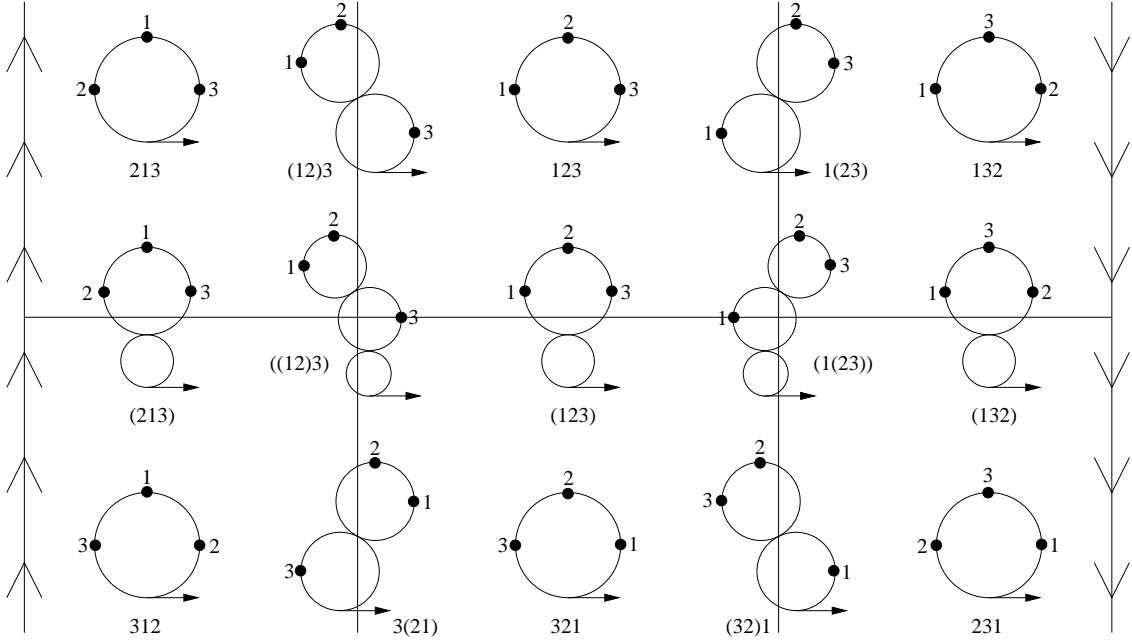,height=8.5cm}
    \caption{The moduli space $\widetilde M_0^4$ is a M\" obius band.}
    \label{fig:modulispace}
  \end{center}
\end{figure}

The entire moduli space has a natural stratification according to the number of double points.  We call a connected component of a stratum a {\bf cell}.
Figure \ref{fig:modulispace} shows a picture of this cell complex\footnote{Unfortunately, this is not a CW-complex, since the closures of the cells are not compact.} for $ n = 3 $.

\subsection{Partial bracketings}
Davis-Januszkiewicz-Scott \cite{DJS}, Devadoss \cite{Dev}, and Kapranov \cite{Kap} have all studied this cell structure on $ \widetilde{M}_0^{n+1} $. Actually, Kapranov studied the associated $\mathbb Z/2$-bundle over $\overline M_0^{n+1}$ (which he denotes $\widetilde {S}^{n-2}$), but his results can easily be adapted to our situation.

An {\bf ordered partial bracketing} on $ n $ letters is an ordering of $ 1 \dots n $ along with an insertion of some (perhaps none) meaningful non-trivial brackets.  For example, $ ((31)(6(54))2) $ is an ordered partial bracketing on $ 6 $ letters.  Two ordered partial bracketings are said to be {\bf equivalent} if one can be attained from the other by reversing the ordering inside brackets.  For example $ 2(14(57)6)(83) $ and $ 2(6(57)41)(38) $ are equivalent.  Let $ F_n $ denote the set of equivalence classes of ordered partial bracketing.

We define a poset structure on $ F_n $ by setting $ \alpha \le \beta $ if there exist representatives $ a $ for $ \alpha $  and $ b $ for $ \beta $ such that $ a $ can be attained from $ b $ by inserting brackets.  Note that there is a natural bijection between $ S_n $ and the  maximum elements of this poset (namely those bracketings with no brackets).  If $ \rho \in S_n $, we will use $ \rho $ for the corresponding maximal element of $ F_n$.  

Each $ \beta \in F_n $ corresponds to a cell of $ \widetilde{M}_0^{n+1} $ consisting of those cacti which have one component for each bracketing and where the order on a component is the order inside the bracket.  For example $ \beta = [2(14(57)6)(83)] $ corresponds to the type of cactus shown in Figure \ref{fig:cactus}. In Figure \ref{fig:modulispace}, we see the bracketing corresponding to (most of) the cells of $\widetilde M_0^4$.

\begin{Theorem}[\cite{Dev}]
Under this bijection between $ F_n $ and the cells of $ \widetilde{M}_0^{n+1} $, the 
 poset structure on $ F_n $ corresponds to the closure poset structure of the cell complex.
\end{Theorem}

The codimension 1 cells correspond to classes of ordered partial bracketings with only one pair of brackets.  For a given top dimensional cell, there is a bijection between its facets and pairs $ (p,q) $ with $ 1 \le p < q \le n$ (corresponding to inserting a bracket to the left of the $ p$-th entry and to the right of the $ q$-th entry).  Moreover, suppose that top dimensional cells $ \gamma$ and $ \gamma' $ are joined by a facet $ (p,q) $.  Then $ \widehat{s}_{p,q} \gamma = \gamma' $ in $S_n$.

\subsection{Admissible paths}
Fix a base point $*$ in the top dimensional cell corresponding to $1\in S_n$.
We say that a path $ P $ in $ \widetilde{M}_0^{n+1} $ is {\bf admissible} if it goes from $\gamma(*)$ to $\gamma'(*)$ and it is transverse to the stratification. In particular, it doesn't intersect the cells of codimension $>1$.

To each admissible path $ P $ we associate an element $ \overline{P} $ of the cactus group $J_n$, by multiplying together appropriate $ s_{p,q} $ according to which codimension 1 cells we cross. If $P$ goes from $\gamma(*)$ to $\gamma'(*)$ then the image of $ \overline{P} $ in $S_n$ is $\gamma'\gamma^{-1}$.

\begin{Theorem} \label{pi1}
Let $ P, Q $ be admissible paths.  Suppose that $ P $ and $ Q $ are homotopic (with fixed endpoints) as curves in $ \widetilde{M}_0^{n+1} $.  Then $ \overline{P} = \overline{Q} $.
In particular, we have an exact sequence of groups
\begin{equation*}
1 \rightarrow \pi_1(\widetilde{M}_0^{n+1}) \rightarrow J_n \rightarrow S_n \rightarrow 1.
\end{equation*}

\end{Theorem}

In analogy with the pure braid group, we call $ \pi_1(\widetilde{M}_0^{n+1})$ the {\bf pure cactus group} and denote it $ PJ_n $.

Using methods of CAT(0) geometry and a result of Gromov, Davis-Januszkiewicz-Scott proved in \cite{DJS} that $\widetilde M_0^{n+1}$ is actually a classifying space for the pure cactus group $PJ_n$. The classifying space of the cactus group can then be identified with the orbifold $[\widetilde M_0^{n+1} / S_n]$ (not to be confused with the space $\widetilde M_0^{n+1} / S_n$).

\begin{proof}
By transversality, if two admissible paths are homotopic, then they are homotopic by a homotopy which is transverse to the stratification. Hence all homotopy relations among admissible paths will arise from passing our paths through codimension 2 cells.  

Two codimension 1 cells $a_1 \cdots (a_p \cdots a_q) \cdots a_n, a_1 \cdots (a_k \cdots a_l) \cdots a_n $ can only have a facet in common if $p<q$ and $k<l$ are disjoint or if one contains the other. 

I
\put(0,-40){\line(1,0){400}}
\put(200,-8){\line(0,-1){60}}
\put(15,-20){$.\,.\, a_k\ldots a_l\ldots a_p\ldots a_q\,.\,.$}
\put(145 ,-20){$.\,.\, a_k\ldots a_l\ldots (a_p\ldots a_q)\,.\,.$}
\put(285  ,-20){$.\,.\, a_k\ldots a_l\ldots a_q\ldots a_p\,.\,.$}
\put(10,-40){$.\,.\, (a_k\ldots a_l)\ldots a_p\ldots a_q\,.\,.$}
\put(140 ,-40){$.\,.\, (a_k\ldots a_l)\ldots (a_p\ldots a_q)\,.\,.$}
\put(280  ,-40){$.\,.\, (a_k\ldots a_l)\ldots a_q\ldots a_p\,.\,.$}
\put(15,-60){$.\,.\, a_l\ldots a_k\ldots a_p\ldots a_q\,.\,.$}
\put(145 ,-60){$.\,.\, a_l\ldots a_k\ldots (a_p\ldots a_q)\,.\,.$}
\put(285  ,-60){$.\,.\, a_l\ldots a_k\ldots a_q\ldots a_p\,.\,.$}
f they are disjoint then we see the local picture:

\noindent which gives the relation $ s_{k,l} s_{p,q} = s_{p,q} s_{k,l} $.  

I
\put(0,-40){\line(1,0){400}}
\put(200,-8){\line(0,-1){60}}
\put(15,-20){$.\,.\, a_p\,.\ldots a_k\ldots a_l\,.\,.\, a_q\,.\,.$}
\put(145 ,-20){$.\,.\, (a_p\,.\ldots a_k\ldots a_l\,.\,.\, a_q)\,.\,.$}
\put(285  ,-20){$.\,.\, a_q\,.\,.\, a_l\ldots a_k\ldots.\, a_p\,.\,.$}
\put(10,-40){$.\,.\, a_p\,.\ldots (a_k\ldots a_l)\,.\,.\, a_q\,.\,.$}
\put(140 ,-40){$.\,.\, (a_p\,.\ldots (a_k\ldots a_l)\,.\,.\, a_q)\,.\,.$}
\put(280  ,-40){$.\,.\, a_q\,.\,.\, (a_k\ldots a_l)\ldots.\, a_p\,.\,.$}
\put(15,-60){$.\,.\, a_p\,.\ldots a_l\ldots a_k\,.\,.\, a_q\,.\,.$}
\put(145 ,-60){$.\,.\, (a_p\,.\ldots a_l\ldots a_k\,.\,.\, a_q)\,.\,.$}
\put(285  ,-60){$.\,.\, a_q\,.\,.\, a_k\ldots a_l\ldots.\, a_p\,.\,.$}
f $ (p, q) $ contains $ (k,l) $ then the local picture is:

\noindent which gives the relation $ s_{p,q} s_{k,l} = s_{m,n} s_{p,q}$, where $ m = \widehat{s}_{p,q}(l)$ and $ n = \widehat{s}_{p,q}(k)$. We see $s_{m,n}$ appear because the numbers $a_k\ldots a_l$ end up in positions $m\ldots n$ after the interval $a_p\ldots a_q$ has been flipped.
This proves that P and Q are homotopic iff $ \overline{P} = \overline{Q} $.

The homomorphism  $J_n \to S_n$ is clearly surjective since the elements $\widehat s_{p,q}$ generate $S_n$.
Representing elements of the cactus group $J_n$ by admissible paths, we see
that the kernel corresponds exactly of those paths which are loops. This implies the result on $\pi_1$.
\end{proof}

Combining Theorems \ref{pi1} and \ref{cactusact}, we see that we can interpret paths in the moduli space as giving us morphisms between multiple tensor products. As we pass through codimension 1 cells, we apply the corresponding reversal to the tensor product.

\begin{Example}
If we consider the path:

\centerline{\epsfig{file=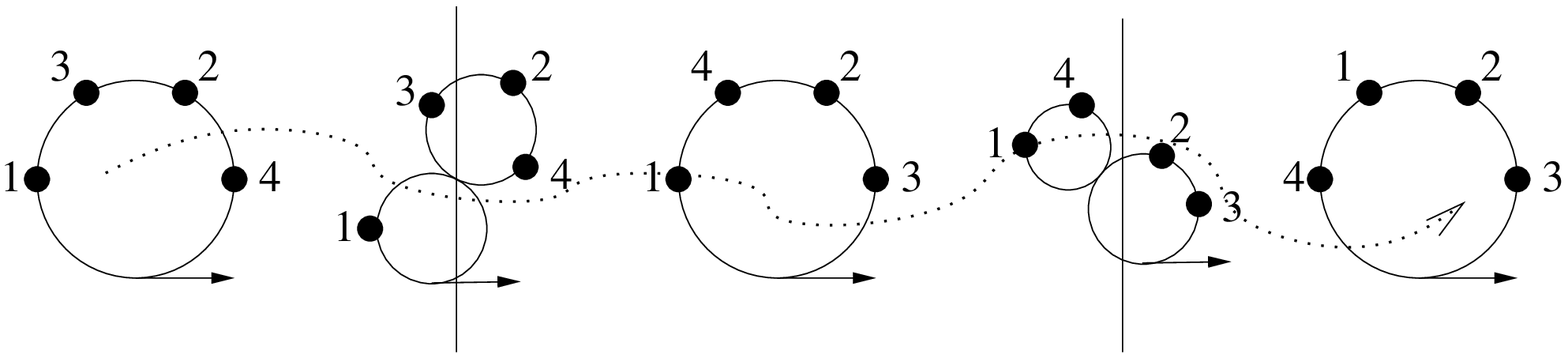,height=2.8cm}}

\noindent we get the morphism:
$$A\otimes C\otimes B\otimes D \xrightarrow{s_{2,4}} A\otimes D\otimes B\otimes C \xrightarrow{s_{1,2}} D\otimes A\otimes B\otimes C.$$
\end{Example}

\subsection{Operads}
The language of operads is useful for restating (and extending) these results.  See \cite{MSS} for the definition of an operad.
First, we note that the spaces $ \widetilde{M}_0^{n+1} $ form an operad. The $i$th composition map 
\begin{equation} \label{eq:opspaces}
\circ_i : \widetilde{M}_0^{n+1} \times \widetilde{M}_0^{k +1} \rightarrow \widetilde{M}_0^{n+k} \\
\end{equation}
takes two cacti $X$ and $Y$ and glues the base point of $Y$ to the $i$th marked point of $X$.
If the base point component of $Y$ had only two special points, it gets collapsed.

Corresponding to this operad in spaces, we get an operad in groupoids by taking the fundamental groupoids.
This construction is entirely analogous to the construction of the braid operad, as noted in \cite{M}.

Let $ G_n = \Pi_1(\widetilde{M}_0^{n+1}, S_n \cdot *) $ be the fundamental groupoid of $\widetilde{M}_0^{n+1}$ 
relative to the basepoints $ S_n \cdot * $.
It is also the action groupoid for the action of $ J_n $ on the symmetric group.
An arrow $(\rho:\pi\to\widehat\rho\pi)\in G_n$ is denoted by $(\rho\,;\pi)$.

The map (\ref{eq:opspaces}) does not send the basepoints of $\widetilde{M}_0^{n+1} \times \widetilde{M}_0^{k+1}$ to 
the basepoints of $\widetilde{M}_0^{n+k}$, so we choose small paths between $(\pi\cdot *)\circ_i(\pi'\cdot *)$ and 
$(\pi\circ_i\pi')\cdot *$. Here, the second $\circ_i$ refers to the usual operad structure on the symmetric groups.
Using these paths, (\ref{eq:opspaces}) induces composition maps
\begin{equation*}
\circ_i : G_n \times G_{k} \rightarrow G_{n+k-1}
\end{equation*}
which endow the groupoids $G_n$ with the structure of an operad.
These composition maps can also be described algebraically on the generators:
\begin{equation*}
\begin{split}
&(1;\pi)\circ_i(s_{p,q};\pi')=(s_{\pi(i)+p-1,\pi(i)+q-1};\pi\circ_i\pi'),\\
&(s_{p,q};\pi)\circ_i(1;\pi')=\begin{cases}
(s_{p,q};\pi\circ_i\pi')&\text{if}\quad \pi(i) > q,\\
(s_{p,q+k-1}s_{\pi(i),\pi(i)+k-1};\pi\circ_i\pi')&\text{if}\quad p\le \pi(i)\le q,\\
(s_{p+k-1,q+k-1};\pi\circ_i\pi')&\text{if}\quad \pi(i) < p.
\end{cases}
\end{split}
\end{equation*}

At this point, one can ask if this operad structure is compatible with the action $\tau$ constructed in Theorem \ref{cactusact}.  
For example, if $(\rho;1)\circ_i(\rho';1)=(\rho'';1)$,
we can consider the following two maps out of the repeated tensor product $ A_1 \cdots A_{n+k-1} $:
\begin{equation*}
\tau \big( \rho ; A_1, \dots, A_{i-1},  \tau(\rho'; A_{i}, \dots, A_{i+k-1}), A_{i+k}, \dots, A_{n+k-1} \big)
\quad\text{and}\quad 
\tau \big(\rho''; A_1, \dots, A_{n+k-1}) \big).
\end{equation*}
Using Lemma \ref{alterdef} and the above description of $ \circ_i $, it is easy to check that these two maps agree.  
 
Thus we obtain the following enhancement of Theorem \ref{cactusact}, which was suggested to us by one of our referees.

\begin{Theorem}
Any coboundary category carries an action (in the weak sense) of the operad $(G_n)$.
\end{Theorem}

\subsection{Further directions}

Analogous to the relation between the braid group and knot theory, the cactus group may be related to some kind of surface topology.  Elements of the pure cactus group can be represented by loops in $\overline{M}_0^{n+1}$.  Given such a loop, we can pullback the tautological curve bundle to get a family of marked stable genus 0 curves over $ S^1 $.  This gives us a surface with $n+1 $ marked lines.  A homotopy of loops gives a cobordism between the corresponding surfaces which preserves the marked lines.  Hence elements of the cactus group correspond to certain surfaces with marked lines modulo certain cobordisms.  It would be interesting to further investigate this connection between the cactus group and 2-dimensional topology.



Drinfel'd's unitarization construction (see section \ref{se:drin}) provides a map $\widehat{PJ_n}\to\widehat{PB_n}$ between the prounipotent completions of the pure cactus group and of the pure braid group. It would be of great interest to learn more about this map (is it injective?).  It would also be of interest to give a geometric version of this map: given a vector bundle with flat unipotent connection on $ \mathbb{C}^n \setminus \Delta  $, we would like to build a vector bundle with flat unipotent connection on $ \widetilde{M}_0^{n+1} $.  In particular, we would like to do this for the Knizhnik-Zamolodchikov connection.  In a similar vein, one can ask to determine the structure of the Lie algebra $\mathrm{Lie}(\widehat{PJ_n})$ analogous to Drinfel'd's computation of $\mathrm{Lie}(\widehat{PB_n})$.

\end{document}